\documentclass[11pt]{amsart}

\newif\ifdraft\draftfalse

\usepackage{amssymb}
\usepackage{amsthm}
\usepackage{eucal}
\usepackage[english]{babel}
\usepackage{nicefrac}
\usepackage{times}
\usepackage{latexsym,amscd,amsmath,amsthm,amssymb}

\makeatletter
\def\@begintheorem#1#2[#3]{%
    \def\naam{#1}
  \deferred@thm@head{\the\thm@headfont \thm@indent
    \@ifempty{#1}{\let\thmname\@gobble}{\let\thmname\@iden}%
    \@ifempty{#2}{\let\thmnumber\@gobble}{\let\thmnumber\@iden}%
    \@ifempty{#3}{\let\thmnote\@gobble}{\let\thmnote\@iden}%
    \thm@swap\swappedhead\thmhead{#1}{#2}{#3}%
    \the\thm@headpunct
    \thmheadnl 
    \hskip\thm@headsep
  }%
  \ignorespaces}
\makeatother

\newcommand{\kantlijndraft}[1]{\ifdraft\hspace{-\lastskip}%
\vadjust{\vspace{-1mm}\smash{\llap{{\tt #1}\hspace{8mm}}}\vspace{1mm}}\fi}

\def\voegToe#1#2#3{\immediate\write1{\string\newlabel{#1}{{#2}{#3}}}}

\newcommand{\thlabel}[1]{\voegToe{#1}{\naam\noexpand~\thetheorem}{\thepage}\kantlijndraft{#1}}

\makeatletter
\renewcommand{\label}[1]{\voegToe{#1}{\@currentlabel}{\thepage}\kantlijndraft{#1}}
\makeatother

\newtheorem{theorem}{Theorem}[section]
\newtheorem{lemma}[theorem]{Lemma}
\newtheorem{corollary}[theorem]{Corollary}
\newtheorem{question}[theorem]{Question}

\newtheorem{proposition}[theorem]{Proposition}
\theoremstyle{definition}
\newtheorem{example}[theorem]{Example}
\newtheorem{definition}[theorem]{Definition}
\theoremstyle{remark}

\numberwithin{equation}{section}

\newtheorem{claim2}{\sc Claim}



\newcommand{\sse}{\subseteq}						
\newcommand{\minus}{\backslash}						
\newcommand{\Un}{\bigcup}							
\newcommand{\un}{\cup}								
\newcommand{\Meet}{\bigcap}							
\newcommand{\meet}{\cap}							
\newcommand{\es}{\varnothing}						

\newcommand{\cl}[1]{\ensuremath{\overline{#1}}}

\newcommand{\scr}[1]{\ensuremath{\mathcal{#1}}}

\def\cprime{$'$}

\def\sapirovskii{{\v{S}}apirovski{\u\i}}
\def\arhangelskii{Arhangel{\cprime}ski{\u\i}}

\def\juhasz{Juh{\'a}sz}

\def\ot {\mathrm{ot}}

\begin{document}

\title{New Bounds on the Cardinality of Hausdorff Spaces and Regular Spaces}

\author{Nathan Carlson}\address{Department of Mathematics, California Lutheran University, 60 W. Olsen Rd, MC 3750, 
Thousand Oaks, CA 91360 USA}
\email{ncarlson@callutheran.edu}

\begin{abstract}
Using weaker versions of the cardinal function $\psi_c(X)$, we derive a series of new bounds for the cardinality of Hausdorff spaces and regular spaces that do not involve $\psi_c(X)$ nor its variants at all. For example, we show if $X$ is regular then $|X|\leq 2^{c(X)^{\pi\chi(X)}}$ and $|X|\leq 2^{c(X)\pi\chi(X)^{ot(X)}}$, where the cardinal function $ot(X)$, introduced by Tkachenko, has the property $ot(X)\leq\min\{t(X),c(X)\}$. It follows from the latter that a regular space with cellularity at most $\mathfrak{c}$ and countable $\pi$-character has cardinality at most $2^\mathfrak{c}$. For a Hausdorff space $X$ we show $|X|\leq 2^{d(X)^{\pi\chi(X)}}$, $|X|\leq d(X)^{\pi\chi(X)^{ot(X)}}$, and $|X|\leq 2^{\pi w(X)^{dot(X)}}$, where $dot(X)\leq\min\{ot(X),\pi\chi(X)\}$. None of these bounds involve $\psi_c(X)$ or $\psi(X)$. By introducing the cardinal functions $w\psi_c(X)$ and $d\psi_c(X)$ with the property $w\psi_c(X)d\psi_c(X)\leq\psi_c(X)$ for a Hausdorff space $X$, we show $|X|\leq\pi\chi(X)^{c(X)w\psi_c(X)}$ if $X$ is regular and $|X|\leq\pi\chi(X)^{c(X)d\psi_c(X)w\psi_c(X)}$ if $X$ is Hausdorff. This improves results of~\sapirovskii~and Sun. It is also shown that if $X$ is Hausdorff then $|X|\leq 2^{d(X)w\psi_c(X)}$, which appears to be new even in the case where $w\psi_c(X)$ is replaced with $\psi_c(X)$. 
Compact examples show that $\psi(X)$ cannot be replaced with $d\psi_c(X)w\psi_c(X)$ in the bound $2^{\psi(X)}$ for the cardinality of a compact Hausdorff space $X$. Likewise, $\psi(X)$ cannot be replaced with $d\psi_c(X)w\psi_c(X)$ in the~\arhangelskii- \sapirovskii~bound $2^{L(X)t(X)\psi(X)}$ for the cardinality of a Hausdorff space $X$. Finally, we make several observations concerning homogeneous spaces in this connection.

\end{abstract}



\maketitle

\section{Introduction.}

Many bounds on the cardinality of a topological space $X$ involve the pseudocharacter cardinal function, $\psi(X)$, or its variant, the closed pseudocharacter $\psi_c(X)$. (See Definition~\ref{psidefinition}). Examples include the~\arhangelskii- \sapirovskii~bound $2^{L(X)t(X)\psi(X)}$ for the cardinality of a Hausdorff space~\cite{Sap1972} and~\sapirovskii's bound $\pi\chi(X)^{c(X)\psi(X)}$ for the cardinality of a regular space~\cite{Sap1974}. This latter result was generalized by Sun~\cite{Sun88} to the class of Hausdorff spaces who showed that $|X|\leq\pi\chi(X)^{c(X)\psi_c(X)}$ if $X$ is Hausdorff. Another example is the Bella-Cammaroto~\cite{BelCam88} bound $d(X)^{t(X)\psi_c(X)}$ for the cardinality of a Hausdorff space, which was improved by the author in~\cite{Car2018} by replacing the tightness $t(X)$ with the weak tightness $wt(X)$ (Definition~\ref{weaktightness}). Other examples include the bounds $\pi w(X)^{ot(X)\psi_c(X)}$ and $\pi\chi(X)^{aL_c(X)ot(X)\psi_c(X)}$ for the cardinality of a Hausdorff space given by Gotchev, Tkachenko, and Tkachuk in~\cite{GTT16}, and the bounds $\pi\chi(X)^{wL(X)ot(X)\psi_c(X)}$ and $2^{wL(X)wt(X)\psi_c(X)}$ for the cardinality of any Hausdorff space with a compact $\pi$-base; that is, a $\pi$-base consisting of elements with compact closure. These bounds were given by Bella, the author, and Gotchev in~\cite{BCG2022b}.

In this study we introduce two cardinal functions for a Hausdorff space $X$ that are weaker than the closed pseudocharacter, $w\psi_c(X)$ and $d\psi_c(X)$ (Definitions~\ref{wpsi} and~\ref{dpsi}), with the property $w\psi_c(X)d\psi_c(X)\leq\psi_c(X)$. We improve the results of~\sapirovskii~and Sun by showing $|X|\leq\pi\chi(X)^{c(X)w\psi_c(X)}$ if $X$ is regular, and $|X|\leq\pi\chi(X)^{c(X)d\psi_c(X)w\psi_c(X)}$ if $X$ is Hausdorff. This is achieved by establishing that $|X|\leq|RO(X)|^{w\psi_c(X)}$, where $RO(X)$ is the collection of regular open subsets of $X$, and, using a modified version of Sun's closing-off argument, that $d(X)\leq\pi\chi(X)^{c(X)d\psi_c(X)}$ for any Hausdorff space $X$. Utilizing upper bounds on $w\psi_c(X)$, it follows that if $X$ is regular then $|X|\leq 2^{c(X)^{\pi\chi(X)}}$ and $|X|\leq 2^{c(X)\pi\chi(X)^{ot(X)}}$. The latter result implies that a regular space with cellularity at most $\mathfrak{c}$, the cardinality of the continuum, and countable $\pi$-character has cardinality at most $2^\mathfrak{c}$. As their proofs rely on sophisticated closing-off arguments, these results are deemed ``difficult" in the sense of Hodel~\cite{Hodel}, section 4. Notice that the last two bounds mentioned do not involve any notion of pseudocharacter.

Using one-to-one map arguments, we can replace $\psi_c(X)$ with $w\psi_c(X)$ in other cardinality bounds. We show $|X|\leq d(X)^{wt(X)w\psi_c(X)}$, $|X|\leq \pi w(X)^{ot(X)w\psi_c(X)}$, and $|X|\leq 2^{d(X)w\psi_c(X)}$ if $X$ is Hausdorff. The latter result appears to be new in the literature even when $w\psi_c(X)$ is replaced with $\psi_c(X)$, although it was known that $|X|\leq 2^{d(X)\psi(X)}$ if $X$ is regular. Using upper bounds on $w\psi_c(X)$, it follows that if $X$ is Hausdorff then $|X|\leq 2^{d(X)^{\pi\chi(X)}}$, $|X|\leq d(X)^{\pi\chi(X)^{ot(X)}}$, and $|X|\leq 2^{\pi w(X)^{dot(X)}}$, where $dot(X)\leq\min\{ot(X),\pi\chi(X)\}$. These curious cardinality bounds again do not involve the pseudocharacter or its variations.

We show in Example~\ref{compactexample} that for every infinite cardinal $\kappa$ there is a compact Hausdorff space $X$ of countable tightness and countable $\pi$-character such that $d\psi_c(X)w\psi_c(X)=\omega$, $\psi(X)\geq\kappa^+$, and $|X|=2^{2^\kappa}$. This shows that $\psi(X)$ cannot be replaced with $d\psi_c(X)w\psi_c(X)$ in the well-known bound $2^{\psi(X)}$ for the cardinality of a compact Hausdorff space. Furthermore, as these examples are countably tight, we see that $\psi(X)$ cannot be replaced with $d\psi_c(X)w\psi_c(X)$ in the~\arhangelskii-~\sapirovskii~bound $2^{L(X)t(X)\psi(X)}$ for the cardinality of a Hausdorff space. In addition, this shows that the bound $\pi\chi(X)^{aL_c(X)ot(X)\psi_c(X)}$~\cite{GTT16} for the cardinality of a Hausdorff space is not valid if $\psi_c(X)$ is replaced with $w\psi_c(X)$. 

In $\S 5$ we make several observations concerning homogeneous Hausdorff spaces and the invariants $d\psi_c(X)$ and $w\psi_c(X)$. First, we observe that $d\psi_c(X)=\psi_c(X)$ for homogeneous Hausdorff spaces. Second, we observe that if $X$ is regular and homogeneous, then $|X|\leq\pi\chi(X)^{c(X)q\psi(X)}$, where $q\psi(X)$ was introduced by Ismail in~\cite{Ism81} with the property $q\psi(X)\leq\min\{\psi_c(X),\pi\chi(X)\}$. (We show further that $q\psi(X)\leq w\psi_c(X)$). We give examples of $\sigma$-compact homogeneous spaces $X$ for which $|X|>2^{w\psi_c(X)}$, demonstrating that $\psi(X)$ cannot be replaced with $w\psi_c(X)$ in the bound $2^{\psi(X)}$ for the cardinality of a $\sigma$-compact space. Nonetheless, we ask the following question: is the cardinality of a compact homogeneous Hausdorff space bounded by $2^{w\psi_c(X)}$? Finally, we note two new bounds for the cardinality of a homogeneous Hausdorff space.

We make no global assumptions of any separation axiom on a topological space in this paper. For definitions not given here, see~\cite{Engelking} and~\cite{Juh80}.

\section{Definitions and preliminary inequalities.}

In this section we give definitions of several cardinal functions used in this paper and prove basic cardinal inequalities that establish interrelationships. We also give definitions of other notions used in this paper.

\begin{definition}\label{psidefinition}
Given a point $x$ in a $T_1$ space $X$, a collection of open sets $\scr{V}$ is a \emph{pseudobase} at $x$ if $\{x\}=\Meet\scr{V}$. If $X$ is additionally Hausdorff, we say $\scr{V}$ is a \emph{closed pseudobase} at $x$ if $\{x\}=\Meet\scr{V}=\Meet_{V\in\scr{V}}\cl{V}$. We define the \emph{pseudocharacter of $x$ in $X$}, denoted by $\psi(x,X)$, to be the least infinite cardinal $\kappa$ such that $x$ has a pseudobase of cardinality $\kappa$. The \emph{closed pseudocharacter of $x$ in $X$}, denoted by $\psi_c(x,X)$, is defined to be the least infinite cardinal $\kappa$ such that $x$ has a closed pseudobase of cardinality $\kappa$. The \emph{pseudocharacter} and \emph{closed pseudocharacter} of the space $X$ are defined by $\psi(X)=\sup\{\psi(x,X):x\in X\}$ and $\psi_c(X)=\sup\{\psi_c(x,X):x\in X\}$, respectively.
\end{definition}

It is clear that $\psi(X)\leq\psi_c(X)\leq\chi(X)$ for any Hausdorff space $X$. Furthermore, if $X$ is regular, it can be seen that $\psi(X)=\psi_c(X)$. 

For a Hausdorff space $X$, we define two new cardinal functions, $w\psi_c(X)$ and $d\psi_c(X)$, each of which are bounded above by $\psi_c(X)$. 

\begin{definition}\label{wpsi}
Let $X$ be a Hausdorff space and let $x\in X$. A collection of open sets $\scr{V}$ is a \emph{weak closed pseudobase at }$x$ if $\{x\}=\Meet_{V\in\scr{V}}\cl{V}$. (Note that it is not necessarily the case that $x\in V$ for any $V\in\scr{V}$). We define $w\psi_c(x,X)$ to be the least infinite cardinal $\kappa$ such that $x$ has a weak closed pseudobase of cardinality $\kappa$. The \emph{weak closed pseudocharacter} $w\psi_c(X)$ is defined as $w\psi_c(X)=\sup\{w\psi_c(x,X):x\in X\}$.
\end{definition}

\begin{definition}\label{dpsi}
Let $X$ be a Hausdorff space. The \emph{dense closed pseudocharacter} $d\psi_c(X)$ of $X$ is defined as the least infinite cardinal $\kappa$ such that $X$ has a dense set $D$ such that $\psi_c(d,X)\leq\kappa$ for every $d\in D$.
\end{definition}

Observe that any space $X$ with a dense set of isolated points has $d\psi_c(X)=\omega$. More generally, any Hausdorff space with a dense set $D$ such that each $d\in D$ has a countable neighborhood base, has $d\psi_c(X)=\omega$. 

A further weakening of the closed pseudocharacter for a Hausdorff space, denoted by $q\psi(X)$, was given by Ismail in~\cite{Ism81}. Ismail showed that if $X$ is Hausdorff then $q\psi(X)\leq\psi_c(X)$ and $q\psi(X)\leq\pi\chi(X)$. We show further in Proposition~\ref{theProp}(6) that $q\psi(X)\leq w\psi_c(X)$ if $X$ is Hausdorff. 

\begin{definition}[\cite{Ism81}]\label{qpsi}
For a point $x$ in a space $X$, a family $\scr{B}$ of nonempty open subsets of $X$ is a $q$\emph{-pseudobase} of $x$ in $X$ if for each $y\in X$ such that $y\neq x$, there is a subfamily $\scr{C}$ of $\scr{B}$ such that $x\in\overline{\Un\scr{C}}$ and $y\notin\overline{\Un\scr{C}}$. We define the \emph{$q$-pseudocharacter} of $x$ in $X$ by $q\psi(x,X)=\min\{|\scr{B}|:\scr{B}\textup{ is a }q\textup{-pseudobase of }x\textup{ in }X\}$ and the $q$-pseudocharacter of $X$ by $q\psi(X)=\sup\{q\psi(x,X):x\in X\}$.
\end{definition}

The cardinal function $ot(X)$ was defined by Tkachenko in~\cite{Tka83} and used by Gotchev, Tkachenko, and Tkachuk in~\cite{GTT16} as well as by Bella, the author, and Gotchev in~\cite{BCG2022b}. The related function $dot(X)$ was defined in~\cite{GTT16}. This function should be considered ``small" in a sense, as $dot(X)\leq\min\{ot(X),\pi\chi(X)\}$ for any space $X$ (Proposition~\ref{theProp}(1)(2)).

\begin{definition}[\cite{Tka83}]
Let $X$ be a space. The \emph{o-tightness} of $X$, denoted by $ot(X)$, is the least infinite cardinal $\kappa$ such that whenever $x\in\cl{\Un\scr{U}}$, for $x\in X$ and an open family $\scr{U}$, there exists $\scr{V}\sse\scr{U}$ such that $|\scr{V}|\leq\kappa$ and $x\in\cl{\Un\scr{V}}$. The \emph{dense o-tightness} of $X$, denoted by $dot(X)$, is the least infinite cardinal $\kappa$ such that whenever $x\in X=\cl{\Un\scr{U}}$, for an open family $\scr{U}$, there exists $\scr{V}\sse\scr{U}$ such that $|\scr{V}|\leq\kappa$ and $x\in\cl{\Un\scr{V}}$.
\end{definition}

The weak tightness $wt(X)$ was defined by author in~\cite{Car2018} and explored further by Bella and the author in~\cite{BC2020}, and by Bella, the author, and Gotchev in~\cite{BCG2022b}. It was shown in~\cite{BCG2022b} that $ot(X)\leq wt(X)$. In~\cite{Car2018} it was shown that if $X$ is Hausdorff then $|X|\leq 2^{L(X)wt(X)\psi(X)}$, improving the well known result of~\arhangelskii~and~\sapirovskii, and in~\cite{BC2020} it was shown that if $X$ is any compact, homogeneous, Hausdorff space then $|X|\leq 2^{wt(X)\pi\chi(X)}$, improving the result of de la Vega~\cite{DLV2006} that such spaces satisfy $|X|\leq 2^{t(X)}$. (Note $\pi\chi(X)\leq t(X)$ for a compactum $X$). We discuss homogeneity more in $\S 5$. 

\begin{definition}[\cite{Car2018}]\label{weaktightness}
The \emph{weak tightness} $wt(X)$ of a space $X$ is defined as the least
infinite cardinal $\kappa$ for which 
there is a cover $\scr{C}$ of $X$ such that $|\scr{C}|\leq
2^\kappa$ and for every $C\in\scr{C}$, $t(C)\leq\kappa$ and 
$X=cl_{2^\kappa}C$.
\end{definition}

\begin{definition}\label{definehomogeneous}
A space $X$ is \emph{homogeneous} if for every pair of points $x$ and $y$ in $X$ there exists a homeomorphism $h:X\to X$ such that $h(x)=y$. 
\end{definition}

The collection of regular open sets plays a prominent role in this work, especially given Proposition~\ref{ROcardBound} below.

\begin{definition}
Let $X$ be a space. A subset of $X$ is \emph{regular closed} if it is of the form $\cl{U}$ where $U$ is open. Let $RC(X)$ denote the set of regular closed sets of $X$. A subset of $X$ is \emph{regular open} if it is of the form $int\cl{U}$ where $U$ is open. Let $RO(X)$ denote the set of regular open sets. 
\end{definition}

It is not hard to see that the complement of a regular closed set is regular open, and therefore the complement of a regular open set is regular closed. Consequently we have $|RO(X)|=|RC(X)|$ for any space $X$.

Proposition~\ref{theProp} below delineates the basic relationships between the cardinal functions defined above. We give full proofs for completeness. In (1), $ot(X)\leq wt(X)$ was originally shown in~\cite{BCG2022b}. (2) was mentioned in~\cite{GTT16}, (3) was mentioned in~\cite{Tka83}, (4), (5), and (6) are new in this paper, and (7) was mentioned in~\cite{Ism81}. For (8), $|RO(X)|\leq\pi w(X)^{c(X)}$ is due to Efimov~\cite{Efimov}, and $|RO(X)|\leq 2^{d(X)}$ is due to de Groot~\cite{DeG1965}. Regarding (5), we point out that a space $X$ needs to be perfectly normal, high on the separation axiom hierarchy, before $w\psi_c(X)=\psi_c(X)=\psi(X)$.

\begin{proposition}\label{theProp}
Let $X$ be a space. Then,
\begin{enumerate}
\item $dot(X)\leq ot(X)\leq wt(X)\leq t(X)$,
\item $dot(X)\leq\pi\chi(X)$,
\item $ot(X)\leq c(X)$,
\item $d\psi_c(X)w\psi_c(X)\leq\psi_c(X)$ if $X$ is Hausdorff,
\item $w\psi_c(X)=\psi(X)$ if $X$ is perfectly normal and $T_1$, 
\item $q\psi(X)\leq w\psi_c(X)$ if $X$ is Hausdorff,
\item $q\psi(X)\leq\pi\chi(X)$ if $X$ is Hausdorff, and
\item $|RO(X)|\leq\min\{\pi w(X)^{c(X)},2^{d(X)}\}$. 
\end{enumerate}
\end{proposition}
\begin{proof}
For (1), $dot(X)\leq ot(X)$ and $wt(X)\leq t(X)$ are clear. To show $ot(X)\leq wt(X)$, let $\kappa=wt(X)$ and let $\scr{C}$ be a cover witnessing that
$wt(X)=\kappa$. Let $x\in\cl{\Un\scr{U}}$, where 
$\scr{U}$ is a family of open sets. There exists $C\in\scr{C}$
such that $x\in C$. Since $C$ is dense in $X$, we have 
$\cl{\Un\scr{U}}=\cl{\Un\scr{U}\meet C}$. Therefore
$x\in\cl{\Un\scr{U}\meet C}\meet C=cl_C(\Un\scr{U}\meet C)$. 
As $t(C)\leq\kappa$, there exists $A\sse\Un\scr{U}\meet C$ such
that $|A|\leq\kappa$ and $x\in cl_C(A)\sse\cl{A}$. 
Thus, there exists $\scr{V}\in[\scr{U}]^\kappa$ such that
$x\in\cl{\Un\scr{V}}$. Therefore $\ot(X)\leq\kappa$.

For (2), let $\kappa=\pi\chi(X)$ and let $x\in X=\cl{\Un\scr{U}}$ for an open family $\scr{U}$. Let $\scr{B}$ be a local $\pi$-base $x$ such that $|\scr{B}|\leq\kappa$. As $\Un\scr{U}$ is dense in $X$, for each $B\in\scr{B}$ there exists $U_B\in\scr{U}$ such that $B\meet U_B\neq\es$. Let $\scr{V}=\{U_B:B\in\scr{B}\}$ and note $|\scr{V}|\leq|\scr{B}|\leq\kappa$. If $W$ is an open set containing $x$ then there exists $B\in\scr{B}$ such that $B\sse W$. Therefore $\es\neq B\meet U_B\sse W\meet U_B\sse W\meet\Un\scr{V}$. This says $x\in\cl{\Un\scr{V}}$ and that $dot(X)\leq\kappa$.

For (3), let $\kappa=c(X)$ and suppose $x\in\cl{\Un\scr{U}}$ for an open family $\scr{U}$. As $c(X)=\kappa$, by using maximal pairwise disjoint open families there exists $\scr{V}\sse\scr{U}$ such that $x\in\cl{\Un\scr{U}}=\cl{\Un\scr{V}}$ and $|\scr{V}|\leq\kappa$. This shows $ot(X)\leq\kappa$.

(4) follows from the fact that every closed pseudobase is a weak closed pseudobase, and that if $\psi_c(X)=\kappa$ then $X$ has a dense subset $D$ such that $\psi_c(d,X)\leq\kappa$ for every $d\in D$, namely $D=X$.

For (5), first recall that in a perfectly normal space every closed set is a $G_\delta$-set. As $w\psi_c(X)\leq\psi_c(X)=\psi(X)$ for any regular space, it remains to show $\psi(X)\leq w\psi_c(X)$. Let $\kappa=w\psi_c(X)$ and let $x\in X$. There exists an open collection $\scr{V}$ such that $\{x\}=\Meet_{V\in\scr{V}}\cl{V}$ and $|\scr{V}|\leq\kappa$. As $X$ is perfectly normal, for each $\scr{V}\in V$ there exists a countable open family $\scr{U}_V$ such that $\cl{V}=\Meet\scr{U}_V$. Then $\{x\}=\Meet_{V\in\scr{V}}\cl{V}=\Meet_{V\in\scr{V}}\Meet\scr{U}_V$. If $\scr{U}=\Un_{V\in\scr{V}}\scr{U}_V$ then $\{x\}=\Meet\scr{U}$ and $|\scr{U}|\leq\kappa\cdot\omega=\kappa$. This says $\psi(X)\leq\kappa$.

For (6), let $\kappa=w\psi_c(X)$. Fix $x\in X$ and let $\scr{V}$ be a weak closed pseudobase at $x$ such that $|\scr{V}|\leq\kappa$. We show $\scr{V}$ serves as a $q$-pseudobase at $x$. If $y\neq x$, then there exists $V\in\scr{V}$ such that $y\in X\minus\cl{V}$. As $x\in\cl{V}$ and $y\in X\minus\cl{V}$ then $\{V\}$ serves as the family ``$\scr{C}$" in the definition of $q$-pseudobase. This shows $\scr{V}$ is a $q$-pseudobase at $x$ and that $q\psi(X)\leq w\psi_c(X)$.

For (7), let $\kappa=\pi\chi(X)$ and fix $x\in X$. Let $\scr{B}$ be a local $\pi$-base at $x$ such that $|\scr{B}|\leq\kappa$. We show $\scr{B}$ is also a $q$-pseudobase at $x$. Let $y\neq x$. As $x$ is Hausdorff, there exists an open set $U$ containing $x$ such that $y\in X\minus\cl{U}$. Let $\scr{C}=\{B\in\scr{B}:B\sse U\}$. Then $x\in\cl{\Un\scr{C}}$ and $y\in X\minus\cl{U}\sse X\minus\cl{\Un\scr{C}}$. This shows $\scr{B}$ is a $q$-pseudobase at $x$ and that $q\psi(X)\leq\kappa$.

For (8), we first show $|RO(X)|\leq\pi w(X)^{c(X)}$. We do this by showing $|RC(X)|\leq\pi w(X)^{c(X)}$. Let $\scr{B}$ be a $\pi$-base for $X$ such that $\scr{B}=\pi w(X)$. For every $R=\cl{U}\in RC(X)$ there exists a collection $\scr{C}_R\sse\scr{B}$ such that $\Un\scr{C}_R\sse U$, $R=\cl{U}=\cl{\Un\scr{C}_R}$ and $|\scr{C}_R|\leq c(X)$. Define a map $\phi:RC(X)\to[\scr{B}]^{\leq c(X)}$ by $\phi(R)=\scr{C}_R$. If $R,Q\in RC(X)$ and $R\neq Q$, then $\scr{C}_R\neq\scr{C}_Q$ for otherwise $R=\cl{\Un\scr{C}_R}=\cl{\Un\scr{C}_Q}=Q$. This shows $\phi$ is one-to-one and therefore $|RO(X)|=|RC(X)|\leq |\scr{B}|^{c(X)}=\pi w(X)^{c(X)}$.

To show $|RO(X)|\leq 2^{d(X)}$, find a dense subset $D$ such that $|D|=d(X)$. We show $RO(X)\sse\{int\cl{A}:A\in\scr{P}(D)\}$. If $U$ is a regular open set then  $U=int\cl{U}$. As $D$ is dense, we have $\cl{U}=\cl{U\meet D}$ and thus $U=int\cl{U}=int(\cl{U\meet D})$. But $A=U\meet D$ is a subset of $D$. This shows $RO(X)\sse\{int\cl{A}:A\in\scr{P}(D)\}$ and thus $|RO(X)|\leq|\scr{P}(D)|=2^{|D|}=2^{d(X)}$.
\end{proof}

Several of the results in this paper that hold for regular spaces in fact hold if the space has a weaker form of regularity known as \emph{quasiregularity}. We define this notion below.

\begin{definition}\label{QR}
A space $X$ is \emph{quasiregular} if every nonempty open set contains a nonempty regular closed set. 
\end{definition}

Clearly any regular space is quasiregular. Any space with a dense set of isolated points is quasiregular. Thus, the space $\kappa\omega$, the Kat\v etov extension of the natural numbers given in Example~\ref{kappaO}, is an example of a nonregular quasiregular Hausdorff space.

\section{Improvements on cardinality bounds of~\sapirovskii~and Sun.}

\sapirovskii~\cite{Sap1974} showed that $d(X)\leq\pi\chi(X)^{c(X)}$ for any regular space $X$. Charlesworth~\cite{Cha77} gave an alternate proof. The author observed in~\cite{Car2007} that this density bound holds if the space is quasiregular (Definition~\ref{QR}), and not necessarily Hausdorff. We give this proof here for completeness. It is a modified version of the proof of 2.37 in~\cite{Juh80}. The reader should compare this proof with the proof of the bound for the density of any Hausdorff space given in Theorem~\ref{densityHausdorff}.

\begin{theorem}[\cite{Sap1974} for regular spaces, \cite{Car2007} for quasiregular spaces]\label{densityQR}
If $X$ is quasiregular then $d(X)\leq\pi\chi(X)^{c(X)}$.
\end{theorem}

\begin{proof}
For $x\in X$ let $\mathcal{B}_x$ be a local $\pi$-base at $x$ such that $|\mathcal{B}_x|\leq\pi\chi(X)$. For $A\sse X$, define $\mathcal{B}_A$ by $\mathcal{B}_A=\Un\{\mathcal{B}_x:x\in A\}$. We now define a map $G:[X]^{\leq c(X)}\rightarrow [X]^{\leq\pi\chi(X)^{c(X)}}$. For $A\in[X]^{\leq c(X)}$, define
$\mathcal{C}_A=\{\scr{U}\in[\mathcal{B}_A]^{\leq c(X)}:X\minus\cl{\Un\scr{U}}\neq\es\}$. Then $\vert\mathcal{C}_A\vert\leq\vert\mathcal{B}_A\vert^{c(X)}\leq (\vert A\vert\cdot\pi\chi(X))^{c(X)}=\pi\chi(X)^{c(X)}$. Now for each $\scr{U}\in\mathcal{C}_A$ pick $p(\scr{U})\in X\minus\cl{\Un\scr{U}}$ and define $G(A)=\{p(\scr{U}):\scr{U}\in\mathcal{C}_A\}\in[X]^{\leq\pi\chi(X)^{c(X)}}$. We now apply \cite[2.24(a)]{Juh80} to obtain a set $A\in [X]^{\leq\pi\chi(X)^{c(X)}}$ such that $G(B)\sse A$ for all $B\in[A]^{\leq c(X)}$. We say that $A$ is closed with respect to $G$.

We claim that $A$ is dense in $X$. Assume the contrary. Then there exists a nonempty open set $W$ such that $W\sse X\minus A$. As $X$ is quasiregular, there exists a nonempty open set $U$ such that $\cl{U}\sse W\sse X\minus A$. 

Now let $\scr{U}$ be a maximal pairwise disjoint family of members of $\scr{B}_A$ disjoint from $\cl{U}$. Suppose there exists $p\in A\minus\cl{\Un\scr{U}}$. Then $p\notin\cl U$, hence $p\in X\minus((\cl{\Un\scr{U})}\un\cl{U})$. There exists $V\in\scr{B}_p\sse\scr{B}_A$ such that $V\sse X\minus((\cl{\Un\scr{U}})\un\cl{U})$. This contradicts the maximality of $\scr{U}$. Hence $A\sse\cl{\Un\scr{U}}$. But $\vert\scr{U}\vert\leq c(X)$ hence we can find a set $H\in [A]^{\leq c(X)}$ such that $\scr{U}\in[\scr{B}_H]^{\leq c(X)}$. Since  $(\Un\scr{U})\meet U=\es$, we have $X\minus\cl{\Un\scr{U}}\neq\es$ and so $\scr{U}\in\scr{C}_H$. Consequently we have $p(\scr{U})\in G(H)\sse A$ as $A$ is closed with respect to $G$. But $p(\scr{U})\in X\minus\cl{\Un\scr{U}}\sse X\minus A$, which is a contradiction. This shows $A$ is dense in $X$. Therefore, $d(X)\leq\vert A\vert\leq\pi\chi(X)^{c(X)}$.

\end{proof}

We modify the proof of Sun~\cite{Sun88} to establish a bound for the density of any Hausdorff space. This bound uses the invariant $d\psi_c(X)$ (see Definition~\ref{dpsi}). 

\begin{theorem}\label{densityHausdorff}
If $X$ is Hausdorff then $d(X)\leq\pi\chi(X)^{c(X)d\psi_c(X)}$.
\end{theorem}

\begin{proof}
Let $\lambda=\pi\chi(X)$ and $\kappa=c(X)d\psi_c(X)$. For all $x\in X$ let $\scr{B}_x$ be a local $\pi$-base for $x$ such that $|\scr{B}_x|\leq\lambda$. Let $D$ be a dense subset of $X$ such that $\psi_c(d,X)\leq\kappa$ for every $d\in D$.

By transfinite induction we construct a non-decreasing chain of $\{A_\alpha:\alpha<\kappa^+\}$ of subsets of $X$ and a sequence of open collections $\{\scr{B}_\alpha:\alpha<\kappa^+\}$ such that the following properties hold for all $\alpha<\kappa^+$:

\begin{enumerate}
\item $|A_\alpha|\leq\lambda^\kappa$,
\item $|\scr{B}_\alpha|\leq\lambda^\kappa$, and 
\item if $\scr{U}=\{\scr{U}_\gamma:\gamma<\kappa\}\in\left[[\scr{B}_\alpha]^{\leq\kappa}\right]^{\leq\kappa}$ and $X\minus\Un_{\gamma<\kappa}\cl{\Un\scr{U}_\gamma}\neq\es$, then $A_\alpha\minus\Un_{\gamma<\kappa}\cl{\Un\scr{U}_\gamma}\neq\es$. 
\end{enumerate}

Pick $p\in X$. Let $A_0=\{p\}$ and $\scr{B}_0=\scr{B}_p$. Let $0<\alpha<\kappa^+$ and assume that $\{A_\beta:\beta<\alpha\}$ have been constructed. Define $\scr{B}_\alpha=\Un\{\scr{B}_x:x\in\Un_{\beta<\alpha}A_\beta\}$. Then $|\scr{B}_\alpha|\leq\lambda\cdot\lambda^\kappa\cdot\kappa=\lambda^\kappa$. For each $\scr{U}=\{\scr{U}_\gamma:\gamma<\kappa\}\in\left[[\scr{B}_\alpha]^{\leq\kappa}\right]^{\leq\kappa}$ such that $X\minus\Un_{\gamma<\kappa}\cl{\Un\scr{U}_\gamma}\neq\es$, pick $x_\scr{U}\in X\minus\Un_{\gamma<\kappa}\cl{\Un\scr{U}_\gamma}$. Define 
$$A_\alpha=\Un_{\beta<\alpha}A_\beta\un\left\{x_\scr{U}:\scr{U}=\{\scr{U}_\gamma:\gamma<\kappa\}\in\left[[\scr{B}_\alpha]^{\leq\kappa}\right]^{\leq\kappa}\textup{ such that }X\minus\Un_{\gamma<\kappa}\cl{\Un\scr{U}_\gamma}\neq\es\right\}.$$
As $\left|\Un_{\beta<\alpha}A_\beta\right|\leq \lambda^\kappa\cdot\kappa=\lambda^\kappa$ and $\left|\left[[\scr{B}_\alpha]^{\leq\kappa}\right]^{\leq\kappa}\right|\leq \left((\lambda^\kappa)^\kappa\right)^\kappa=\lambda^\kappa$, we see that $|A_\alpha|\leq\lambda^\kappa$. By the way we have constructed $A_\alpha$ we see that (3) is satisfied.

Let $A=\Un_{\alpha<\kappa^+}A_\alpha$. Then $|A|\leq\kappa^+\lambda^\kappa=\lambda^\kappa$. We show that $A$ is dense in $X$. Suppose by way of contradiction that there exists a nonempty open set $U$ such that $U\meet A=\es$. As $D$ is dense in $X$, there exists $d\in U\meet D\sse D\minus A$. As $d\in D$, there is a closed pseudobase $\scr{V}=\{V_\alpha:\alpha<\kappa\}$ at $d$ such that $|\scr{V}|=\kappa$. Then $\{d\}=\Meet\scr{V}=\Meet_{\alpha<\kappa}\cl{V_\alpha}$.

For every $\alpha<\kappa$, let $W_\alpha=X\minus\cl{V_\alpha}$. Then $\cl{W_\alpha}=cl(X\minus\cl{V_\alpha})=X\minus int(\cl{V_\alpha})\sse X\minus V_\alpha\sse X\minus\{d\}$ and so $d\notin\cl{W_\alpha}$. Furthermore, $A\sse X\minus\{d\}=\Un_{\alpha<\kappa}W_\alpha$. For each $\alpha<\kappa$, define $\scr{S}_\alpha=\{B\in\scr{B}_x: x\in W_\alpha\meet A, B\sse W_\alpha\}$. Note $\Un\scr{S}_\alpha\sse W_\alpha$.

We show $A\meet W_\alpha\sse\cl{\Un\scr{S}_\alpha}$ for all $\alpha<\kappa$. Let $x\in A\meet W_\alpha$ and let $T$ be an open set containing $x$. There exists $B\in\scr{B}_x$ such that $B\sse W_\alpha\meet T$. Then $B\in\scr{S}_\alpha$ and $\es\neq B\sse T\meet\Un\scr{S}_\alpha$. This shows $A\meet W_\alpha\sse\cl{\Un\scr{S}_\alpha}$ for all $\alpha<\kappa$.

As $c(X)\leq\kappa$, for each $\alpha<\kappa$ there exists $\scr{U}_\alpha\sse\scr{S}_\alpha$ such that $\cl{\Un\scr{U}_\alpha}=\cl{\Un\scr{S}_\alpha}$ and $|\scr{U}_\alpha|\leq\kappa$. Now, for each $\alpha<\kappa$ note that $\cl{\Un\scr{U}_\alpha}\sse\cl{W_\alpha}\sse X\minus\{d\}$. Therefore, $d\in X\minus\Un_{\alpha<\kappa}\cl{\Un\scr{U}_\alpha}$. Since $\left|\Un_{\alpha<\kappa}\scr{U}_\alpha\right|\leq\kappa\cdot\kappa=\kappa<\kappa^+$, there exists $\delta<\kappa^+$ such that $\scr{U}=\{\scr{U}_\alpha:\alpha<\kappa\}\in\left[[\scr{B}_\delta]^{\leq\kappa}\right]^{\leq\kappa}$.

By (3) above, we have that $x_\scr{U}\in A_{\delta+1}\minus\Un_{\alpha<\kappa}\cl{\Un\scr{U}_\alpha}\sse A\minus\Un_{\alpha<\kappa}\cl{\Un\scr{U}_\alpha}$. This contradicts the fact that $A\sse\Un_{\alpha<\kappa}\cl{\Un\scr{S}_\alpha}=\Un_{\alpha<\kappa}\cl{\Un\scr{U}_\alpha}$. Therefore $A$ is dense in $X$ and
$d(X)\leq |A|\leq\lambda^\kappa=\pi\chi(X)^{c(X)d\psi_c(X)}$.
\end{proof}

One should compare the proofs of Theorems~\ref{densityQR} and~\ref{densityHausdorff}, as they both give bounds for the density of a space. The proofs are sophisticated closing-off arguments, but different. We ask the following:

\begin{question}
Is there a common proof of Theorems~\ref{densityQR} and~\ref{densityHausdorff}?
\end{question}

By~\ref{densityQR} and~\ref{densityHausdorff} we have a class of spaces for which $d(X)\leq\pi\chi(X)^{c(X)}$. This is given in the next corollary.

\begin{corollary}
Let $X$ be a space. If $X$ is either quasiregular, or is Hausdorff with a dense set of points with countable closed pseudocharacter, then $d(X)\leq\pi\chi(X)^{c(X)}$.
\end{corollary}

\begin{example}
In~\cite{CR2008} the author and Ridderbos constructed an involved example under the axiom $\mathfrak{c}^+=2^\mathfrak{c}$ of a c.c.c Urysohn space $Z$ with $\pi$-character $\mathfrak{c}$ and $d(Z)=\mathfrak{c}^+=2^\mathfrak{c}$. Therefore we have $d(Z)=2^\mathfrak{c}>\mathfrak{c}=\pi\chi(Z)^{c(Z)}$. In other words, the cardinal inequality $d(X)\leq\pi\chi(X)^{c(X)}$ is not valid for all Urysohn spaces despite being valid for all quasiregular spaces (Theorem~\ref{densityQR}). Nonetheless, since $Z$ is Hausdorff, by Theorem~\ref{densityHausdorff} we have $\mathfrak{c}^+=d(Z)\leq\pi\chi(Z)^{c(Z)d\psi_c(X)}=\mathfrak{c}^{d\psi_c(Z)}$. This implies $d\psi_c(Z)$ is large, at least uncountable.
\end{example}

The following was proved by Efimov~\cite{Efimov} for regular spaces. We show in fact it works for any quasiregular space (not necessarily Hausdorff).

\begin{corollary}[Efimov~\cite{Efimov} in the regular case]\label{QRRObound}
If $X$ is quasiregular then $|RO(X)|\leq\pi\chi(X)^{c(X)}$.
\end{corollary}

\begin{proof}
Applying Theorem~\ref{densityQR} and the fact that $\pi w(X)=d(X)\pi\chi(X)$, we have $\pi w(X)\leq\pi\chi(X)^{c(X)}$. By Proposition~\ref{theProp}(8), we have
$$|RO(X)|\leq\pi w(X)^{c(X)}\leq\left(\pi\chi(X)^{c(X)}\right)^{c(X)}=\pi\chi(X)^{c(X)}.$$
\end{proof}

In a similar manner, we derive a bound for the cardinality of $RO(X)$ for Hausdorff spaces. 

\begin{corollary}\label{HausRO}
If $X$ is Hausdorff then $|RO(X)|\leq\pi\chi(X)^{c(X)d\psi_c(X)}$.
\end{corollary}

\begin{proof}
Appling Theorem~\ref{densityHausdorff} and the fact $\pi w(X)=d(X)\pi\chi(X)$ for any space, we see that $\pi w(X)\leq\pi\chi(X)^{c(X)d\psi_c(X)}$. By Proposition~\ref{theProp}(8), we have
$$|RO(X)|\leq\pi w(X)^{c(X)}\leq\left(\pi\chi(X)^{c(X)d\psi_c(X)}\right)^{c(X)}=\pi\chi(X)^{c(X)d\psi_c(X)}.$$
\end{proof}

In the following Proposition we show the cardinality of a Hausdorff space $X$ is related to $|RO(X)|$. This simple result is at the core of several of the cardinality bounds in the next section. It appears not to be mentioned in the literature, even in the case where $w\psi_c(X)$ is replaced with the usual $\psi_c(X)$.

\begin{proposition}\label{ROcardBound}
If $X$ is Hausdorff, then $|X|\leq |RO(X)|^{w\psi_c(X)}$. 
\end{proposition}

\begin{proof}
Let $\kappa=w\psi_c(X)$. For each $x\in X$, let $\scr{V}_x$ be a family of open sets such that $\{x\}=\Meet_{V\in\scr{V}_x}\cl{V}$ and $|\scr{V}_x|\leq\kappa$. Define $\phi:X\to[RC(X)]^{\leq\kappa}$ by $\phi(x)=\{\cl{V}:V\in\scr{V}_x\}$. Then $\phi$ is one-to-one: suppose $x\neq y\in X$. Then $\Meet_{V\in\scr{V}_x}\cl{V}=\{x\}\neq\{y\}=\Meet_{V\in\scr{V}_y}\cl{V}$, and thus $\{\cl{V}:V\in\scr{V}_x\}\neq\{\cl{V}:V\in\scr{V}_y\}$. This shows $\phi(x)\neq\phi(y)$ and that $\phi$ is one-to-one. Therefore, $|X|\leq |RC(X)|^\kappa=|RO(X)|^\kappa=|RO(X)|^{w\psi_c(X)}$.
\end{proof}

In 1974 \sapirovskii~\cite{Sap1974} showed that the cardinality of any regular Hausdorff space $X$ is at most $\pi\chi(X)^{c(X)\psi(X)}$. We improve upon this result in the next theorem by showing such spaces have cardinality bounded by $\pi\chi(X)^{c(X)w\psi_c(X)}$. This is a logical improvement for regular spaces as in that case $w\psi_c(X)\leq\psi_c(X)=\psi(X)$. This appears to be the first known improvement on this result of \sapirovskii. In fact, our inequality works for any quasiregular Hausdorff space.

\begin{theorem}\label{ShaImprove}
If $X$ is quasiregular and Hausdorff then $|X|\leq\pi\chi(X)^{c(X)w\psi_c(X)}$.
\end{theorem}

\begin{proof}
By Corollary~\ref{QRRObound} and Proposition~\ref{ROcardBound}, we have
$$|X|\leq |RO(X)|^{w\psi_c(X)}\leq\left(\pi\chi(X)^{c(X)}\right)^{w\psi_c(X)}=\pi\chi(X)^{c(X)w\psi_c(X)}.$$
\end{proof}

In 1988~\cite{Sun88} Sun extended~\sapirovskii's result for regular spaces by showing the cardinality of any Hausdorff space $X$ is bounded by $\pi\chi(X)^{c(X)\psi_c(X)}$. This improved the Hajnal-\juhasz~theorem that such spaces have cardinality bounded by $2^{c(X)\chi(X)}$. We give a logical improvement of Sun's bound below, where $\psi_c(X)$ is replaced with $d\psi_c(X)w\psi_c(X)$, recalling that $d\psi_c(X)w\psi_c(X)\leq\psi_c(X)$ for any Hausdorff space by Proposition~\ref{theProp}(4). Theorem~\ref{SunImprove} appears to be the first improvement of Sun's bound.

\begin{theorem}\label{SunImprove}
If $X$ is Hausdorff then $|X|\leq\pi\chi(X)^{c(X)d\psi_c(X)w\psi_c(X)}$.
\end{theorem}

\begin{proof}
By Proposition~\ref{ROcardBound} and Corollary~\ref{HausRO}, we have
$$|X|\leq |RO(X)|^{w\psi_c(X)}\leq\left(\pi\chi(X)^{c(X)d\psi_c(X)}\right)^{w\psi_c(X)}=\pi\chi(X)^{c(X)d\psi_c(X)w\psi_c(X)}.$$
\end{proof}

By Theorem~\ref{SunImprove}, it follows that any c.c.c Hausdorff space $X$ with $\pi$-character at most $\mathfrak{c}$, a dense set $D$ such that $\psi_c(d,X)\leq\omega$ for all $d\in D$, and a countable weak closed pseudobase at every point, has $|X|\leq\mathfrak{c}$.

The next example, simple but illustrative, shows that for every infinite cardinal $\kappa$ there is a compact Hausdorff space $X$ of countable tightness and countable $\pi$-character such that $d\psi_c(X)w\psi_c(X)=\omega$, $\psi(X)=\psi_c(X)\geq\kappa^+$, and $|X|=2^{2^\kappa}$.

\begin{example}\label{compactexample}
Let $D$ be any infinite discrete space and let $X=D\un\{p\}$ be the one-point compactification of $D$ where $p$ is the point at infinity. We aim to show that $d\psi_c(X)w\psi_c(X)=\omega$. First, note that $d\psi_c(X)=\omega$ as $X$ has a dense set of isolated points. Now, notice that $p$ is in the closure of any infinite subset of $D$ as neighborhoods of $p$ contain all but finitely many elements of $D$. This makes $X$ countably tight and by~\sapirovskii's result that $\pi\chi(Z)\leq t(Z)$ for any compact Hausdorff space $Z$, we have that $X$ has countable $\pi$-character. 

Let $A$ be any countably infinite subset of $D$. For each $x\in A$, let $A_x=A\minus\{x\}$. Then $A_x$ is infinite for each $x\in A$ and thus $p\in\cl{A_x}$ for every $x\in A$. Note further that each $A_x$ is open as it consists of isolated points. Now, every $y\in D\minus A$ is not in $\cl{A_x}$ for any $x\in A$ as $y$ is isolated. Additionally, for each $x\in A$ we have $x\notin\cl{A_x}$. This shows $\{p\}=\Meet_{x\in A}\cl{A_x}$ and as $A$ is countable we see that $\{A_x:x\in A\}$ is a countable weak closed pseudobase at $p$. As $x$ is isolated for every $x\in X\minus\{p\}$, we see that $w\psi_c(X)=\omega$.

Now let $D$ be the discrete space of cardinality $2^{2^\kappa}$ for an infinite cardinal $\kappa$. Then the space $X$ has cardinality $2^{2^\kappa}$ and, since $X$ is compact, we have $|X|\leq 2^{\psi(X)}$. Therefore $2^{2^\kappa}\leq 2^{\psi(X)}$ and $\psi(X)\geq\kappa^+$, for otherwise if $\psi(X)\leq\kappa$ then $2^{2^\kappa}\leq 2^\kappa$, a contradiction. Thus, by the above, for every infinite cardinal $\kappa$ there is a compact space $X$ of countable tightness and countable $\pi$-character such that $d\psi_c(X)w\psi_c(X)=\omega$ and $\psi(X)=\psi_c(X)\geq\kappa^+$. As $|X|=2^{2^\kappa}$, we have that $2^{d\psi_c(X)w\psi_c(X)}$ is not a bound for the cardinality of all compact Hausdorff spaces. Also, since the space $X$ is countably tight, this implies that $2^{L(X)t(X)d\psi_c(X)w\psi_c(X)}$ is not a bound for the cardinality of all Hausdorff spaces, i.e. $\psi(X)$ cannot be replaced with $d\psi_c(X)w\psi_c(X)$ in the~\arhangelskii-~\sapirovskii~bound $2^{L(X)t(X)\psi(X)}$ for Hausdorff spaces. Moreover, it was shown in~\cite{GTT16} that the cardinality of a Hausdorff space is always bounded by $\pi\chi(X)^{aL_c(X)ot(X)\psi_c(X)}$, where $aL_c(X)$, the almost Lindel\"of degree with respect to closed sets, satisfies $aL_c(X)\leq L(X)$. Our compact space $X$ is also a counterexample to replacing $\psi_c(X)$ with $d\psi_c(X)w\psi_c(X)$ in this bound.

A further observation concerns Theorem~\ref{SunImprove}. Applying that theorem to the space $X$, we see $2^{2^\kappa}=|X|\leq\pi\chi(X)^{c(X)d\psi_c(X)w\psi_c(X)}=\omega^{c(X)\cdot\omega}=2^{c(X)}$, implying the cellularity of $X$ must be large. Clearly, it is large since $X$ contains $2^{2^\kappa}$ isolated points. Therefore in Sun's bound $\pi\chi(X)^{c(X)\psi_c(X)}$ for the cardinality of a Hausdorff space, the closed pseudocharacter is larger than necessary. Indeed that cardinal function can be replaced with $d\psi_c(X)w\psi_c(X)$ which, as we see in this example, can be countable. This necessarily forces the cellularity to be large.

\end{example}

\begin{example}\label{kappaO}
Consider $X=\kappa\omega$, the Kat\v etov extension of the natural numbers. Observe that $X$ has a countable dense set of isolated points, and so $X$ is c.c.c, has countable $\pi$-character, and $d\psi_c(X)=\omega$. Note also $|X|=2^\mathfrak{c}$, $\psi_c(X)=\mathfrak{c}$, and that $X$ is Hausdorff. By Theorem~\ref{SunImprove}, $2^\mathfrak{c}=|X|\leq\pi\chi(X)^{c(X)d\psi_c(X)w\psi_c(X)}=2^{w\psi_c(X)}$, and so $w\psi_c(X)$ must be uncountable. As $X$ has a dense set of isolated points, it is quasiregular, implying we could also have arrived at the same conclusion using Theorem~\ref{ShaImprove}.
\end{example}

In~\cite{BCG2022b}, Bella, the author, and Gotchev studied Hausdorff spaces with a \emph{compact $\pi$-base}, that is, a $\pi$-base with elements with compact closure. This class of spaces generalizes both the class of locally compact spaces and the class of spaces with a dense set of isolated points. It was shown in~\cite{BCG2022b} that such spaces are always quasiregular Baire spaces. Therefore, by Theorem~\ref{ShaImprove}, we have the following corollary.

\begin{corollary}\label{cptpibase}
If $X$ is a Hausdorff space with a compact $\pi$-base, then $|X|\leq\pi\chi(X)^{c(X)w\psi_c(X)}$.
\end{corollary}

It was shown in~\cite{BCG2022b} that if $X$ is Hausdorff with a compact $\pi$-base then $|X|\leq 2^{wL(X)wt(X)\psi_c(X)}$ and $|X|\leq\pi\chi(X)^{wL(X)ot(X)\psi_c(X)}$, where $wL(X)$ is the weak Lindel\"of degree of $X$. One may ask if $\psi_c(X)$ can be replaced with $d\psi_c(X)w\psi_c(X)$ in either of these two bounds. Example~\ref{compactexample} provides a negative answer. The example $X$ in~\ref{compactexample} is a compact space (hence it has a compact $\pi$-base) where all the relevant cardinal functions involved in these bounds are countable but we can make $|X|$ as large as necessary. Therefore neither $2^{wL(X)wt(X)d\psi_c(X)w\psi_c(X)}$ nor $\pi\chi(X)^{wL(X)ot(X)d\psi_c(X)w\psi_c(X)}$ are bounds for the cardinality of all Hausdorff spaces with a compact $\pi$-base.

\section{New cardinal inequalities using upper bounds on $w\psi_c(X)$.}

We now obtain several upper bounds for $w\psi_c(X)$ for Hausdorff spaces. These are given in Propositions~\ref{wpc1},~\ref{wpc2},~\ref{wpc3}, and~\ref{wpc4}. These cardinal inequalities are used to generate new bounds for the cardinality of Hausdorff spaces and quasiregular Hausdorff spaces that do not involve the pseudocharacter of a space nor its variants.

\begin{proposition}\label{wpc1}
If $X$ is a Hausdorff space then $w\psi_c(X)\leq\pi\chi(X)^{ot(X)}$.
\end{proposition}

\begin{proof}
Let $x\in X$ and let $\scr{B}$ be a local $\pi$-base at $x$ such that $|\scr{B}|\leq\pi\chi(X)$. As $X$ is Hausdorff, for all $y\neq x$ there exists an open set $U_y$ containing $x$ such that $y\in X\minus\cl{U_y}$. Then $x\in cl\left(\Un\{B\in\scr{B}:B\sse U_y\}\right)\sse\cl{U_y}\sse X\minus\{y\}$. There exists $\scr{U}_y\sse\{B\in\scr{B}:B\sse U_y\}$ such that $|\scr{U}_y|\leq ot(X)$ and $x\in\cl{\Un\scr{U}_y}\sse\cl{U_y}\sse X\minus\{y\}$. This shows $\{x\}=\Meet_{y\neq x}\cl{\Un\scr{U}_y}$ and that $\scr{V}=\{\Un\scr{U}_y:y\neq x\}$ is a weak closed pseudobase at $x$. Now, $\scr{V}\sse\{\Un\scr{C}:\scr{C}\in[\scr{B}]^{\leq\ot(X)}\}$ which implies $|\scr{V}|\leq\left|[\scr{B}]^{\leq\ot(X)}\right|\leq|\scr{B}|^{ot(X)}\leq\pi\chi(X)^{ot(X)}$.
\end{proof}

The next proposition demonstrates that $dot(X)$ ``acts like" $ot(X)$ in a particular situation. We will use this proposition to establish another bound for $w\psi_c(X)$ given in Proposition~\ref{wpc2}.

\begin{proposition}\label{dot}
Let $X$ be a space, let $dot(X)\leq\kappa$, and let $x\in V$ where $V$ is open. If $x\in\cl{\Un\scr{U}}$ where $\scr{U}$ is an family of open sets such that $\cl{\Un\scr{U}}=\cl{V}$, then there exists $\scr{V}\in[\scr{U}]^{\leq\kappa}$ such that $x\in\cl{\Un\scr{V}}$. 
\end{proposition}

\begin{proof}
Let $W=X\minus\cl{V}=X\minus\cl{\Un\scr{U}}$. Then $X=W\un\cl{\Un\scr{U}}\sse\cl{W}\un\cl{\Un\scr{U}}=cl(W\un\Un\scr{U})$. Then $x\in cl(W\un\Un\scr{U})$ and since $dot(X)\leq\kappa$ there exists $\scr{V}\in[\scr{U}]^{\leq\kappa}$ such that $x\in cl(W\un\Un\scr{V})$. As $x\in V$, we have $x\in cl(V\meet(W\un\Un\scr{V}))$. But $V\meet W=\es$ and so $x\in cl(V\meet\Un\scr{V})\sse\cl{\Un\scr{V}}$.
\end{proof}

\begin{proposition}\label{wpc2}
If $X$ is a Hausdorff space then $w\psi_c(X)\leq\pi w(X)^{dot(X)}$.
\end{proposition}

\begin{proof}
We proceed with an argument similar to that in Proposition~\ref{wpc1} and use Proposition~\ref{dot}. Let $\scr{B}$ be a $\pi$-base for $X$ such that $|\scr{B}|=\pi w(X)$ and fix $x\in X$. As $X$ is Hausdorff, for all $y\neq x$ there exists an open set $U_y$ containing $x$ such that $y\in X\minus\cl{U_y}$. Then $x\in cl\left(\Un\{B\in\scr{B}:B\sse U_y\}\right)=\cl{U}_y\sse X\minus\{y\}$. Notice that $cl\left(\Un\{B\in\scr{B}:B\sse U_y\}\right)=\cl{U}_y$ because $\scr{B}$ is a $\pi$-base for $X$. As $x\in U_y$, by Proposition~\ref{dot} there exists $\scr{U}_y\sse\{B\in\scr{B}:B\sse U_y\}$ such that $|\scr{U}_y|\leq dot(X)$ and $x\in\cl{\Un\scr{U}_y}\sse\cl{U_y}\sse X\minus\{y\}$. This shows $\{x\}=\Meet_{y\neq x}\cl{\Un\scr{U}_y}$ and that $\scr{V}=\{\Un\scr{U}_y:y\neq x\}$ is a weak closed pseudobase at $x$. Now, $\scr{V}\sse\{\Un\scr{C}:\scr{C}\in[\scr{B}]^{\leq dot(X)}\}$ which implies $|\scr{V}|\leq\left|[\scr{B}]^{\leq dot(X)}\right|\leq|\scr{B}|^{dot(X)}\leq\pi w(X)^{dot(X)}$.
\end{proof}

\begin{proposition}\label{wpc3}
If $X$ is a Hausdorff space then $w\psi_c(X)\leq 2^{d(X)}$.
\end{proposition}

\begin{proof}
Let $D$ be a dense set such that $|D|=d(X)$ and fix $x\in X$. As $X$ is Hausdorff for all $y\neq x$ there exists an open set $U_y$ containing $x$ such that $y\in X\minus\cl{U_y}$. Then, $\{x\}=\Meet_{y\neq x}\cl{U_y}=\Meet_{y\neq x}\cl{U_y\meet D}$ as $D$ is dense. Then $\scr{D}=\{U_y\meet D:y\neq x\}\sse\scr{P}(D)$ and so $|\scr{D}|\leq |\scr{P}(D)|=2^{|D|}=2^{d(X)}$. 

For all $A\in\scr{D}$ there exists $y_A\neq x$ such that $A=U_{y_A}\meet D$. Also note $x\in\cl{A}$ for all $A\in\scr{D}$. Then, $\{x\}=\Meet_{A\in\scr{D}}\cl{A}=\Meet_{A\in\scr{D}}\cl{U_{y_A}\meet D}=\Meet_{A\in\scr{D}}\cl{U_{y_A}}$. This shows $\scr{V}=\{U_{y_A}:A\in\scr{D}\}$ is a weak closed pseudobase at $X$. Moreover, $|\scr{V}|\leq|\scr{D}|\leq 2^{d(X)}$. This shows $w\psi_c(x,X)\leq 2^{d(X)}$ for any $x\in X$ and that $w\psi_c(X)\leq 2^{d(X)}$. 
\end{proof}

\begin{proposition}\label{wpc4}
If $X$ is a Hausdorff space then $w\psi_c(X)\leq 2^{q\psi(X)}$.
\end{proposition}

\begin{proof}
Let $\kappa=q\psi(X)$ and fix $x\in X$. Let $\scr{B}$ be a $q$-pseudobase at $x$ such that $|\scr{B}|\leq\kappa$. For all $y\neq x$ there exists $\scr{C}_y\in\scr{P}(\scr{B})$ such that $x\in\cl{\Un\scr{C}_y}$ and $y\notin\cl{\Un\scr{C}_y}$. Then $\{x\}=\Meet_{y\neq x}\cl{\Un\scr{C}_y}$ and $\scr{C}=\{\Un\scr{C}_y:y\neq x\}$ is a weak closed pseudobase for $x$. Now, $\scr{C}\sse\{\Un\scr{D}:\scr{D}\in\scr{P}(\scr{B})\}$ and so $|\scr{C}|\leq|\scr{P}(\scr{B})|=2^{|\scr{B}|}=2^\kappa$. This says $w\psi_c(x,X)\leq 2^\kappa$ for all $x\in X$, and so $w\psi_c(X)\leq 2^\kappa$.
\end{proof}

Corollaries~\ref{QRBound} and~\ref{HausBound} below depend on Theorems~\ref{ShaImprove} and~\ref{SunImprove}, which in turn depend on the closing-off arguments in the proofs of Theorems~\ref{densityQR} and~\ref{densityHausdorff}. As the proofs of~\ref{QRBound} and~\ref{HausBound} ultimately depend on these sophisticated techniques, they are ``difficult'' bounds in the sense of Hodel~\cite{Hodel}, section 4. The bounds in Corollary~\ref{QRBound} are ones that do not involve any notion of pseudocharacter.

\begin{corollary}\label{QRBound}
Let $X$ be a quasiregular Hausdorff space. Then, 
\begin{itemize}
\item[(1)] $|X|\leq 2^{c(X)\pi\chi(X)^{ot(X)}}$, 
\item[(2)] $|X|\leq\pi\chi(X)^{c(X)^{q\psi(X)}}$, and 
\item[(3)] $|X|\leq 2^{c(X)^{\pi\chi(X)}}$.
\end{itemize}
\end{corollary}

\begin{proof}
For (1), by Theorem~\ref{ShaImprove} and Proposition~\ref{wpc1}, we have
$$|X|\leq\pi\chi(X)^{c(X)w\psi_c(X)}\leq\pi\chi(X)^{c(X)\pi\chi(X)^{ot(X)}}=2^{c(X)\pi\chi(X)^{ot(X)}}.$$
For (2), by Theorem~\ref{ShaImprove} and Proposition~\ref{wpc4}, we have
$$|X|\leq\pi\chi(X)^{c(X)w\psi_c(X)}\leq\pi\chi(X)^{c(X)\cdot 2^{q\psi(X)}}=\pi\chi(X)^{c(X)^{q\psi(X)}}.$$
(3) follows from (2) and the fact that if $X$ is Hausdorff then $q\psi(X)\leq\pi\chi(X)$ (Proposition~\ref{theProp}(7)).
\end{proof}

It follows from Corollary~\ref{QRBound}(3) that a quasiregular Hausdorff space with cellularity at most $\mathfrak{c}$ and countable $\pi$-character has cardinality at most $2^\mathfrak{c}$. 

\begin{corollary}\label{HausBound}
If $X$ is Hausdorff then $|X|\leq 2^{c(X)d\psi_c(X)\pi\chi(X)^{ot(X)}}$.
\end{corollary}

\begin{proof}
Apply Theorem~\ref{SunImprove} and Proposition~\ref{wpc1}.
\end{proof}

Compare Corollary~\ref{HausBound} with Corollary~\ref{QRBound}(1) for quasiregular Hausdorff spaces. The difference lies in the addition of the function $d\psi_c(X)$ in Corollary~\ref{HausBound}. Recall $d\psi_c(X)$ is countable, for example, when a space $X$ has a dense subspace $D$ with countable neighborhood bases at each point in $D$.

In light of Corollaries~\ref{QRBound}(1) and~\ref{HausBound}, we ask the following.

\begin{question}
If $X$ is Hausdorff, is $|X|\leq 2^{c(X)\pi\chi(X)^{ot(X)}}$?
\end{question}

The bounds in the rest of this section do not rest upon any sophisticated closing-off arguments but rather more straightforward arguments involving one-to-one maps. These are not ``difficult'' bounds in the sense of Hodel~\cite{Hodel}, section 4. Theorem~\ref{dwpsi} below appears to be new in the literature, even replacing $w\psi_c(X)$ with $\psi_c(X)$.

\begin{theorem}\label{dwpsi}
If $X$ is Hausdorff then $|X|\leq 2^{d(X)w\psi_c(X)}$.
\end{theorem}

\begin{proof}
By Propositions~\ref{theProp}(8) and~\ref{ROcardBound}, we have
$$|X|\leq|RO(X)|^{w\psi_c(X)}\leq{\left(2^{d(X)}\right)^{w\psi_c(X)}}=2^{d(X)w\psi_c(X)}.$$
\end{proof}

It follows that the cardinality of any regular Hausdorff space $X$ is at most $2^{d(X)w\psi_c(X)}\leq 2^{d(X)\psi_c(X)}=2^{d(X)\psi(X)}$, a well-known result. (See, for example, Problem 3.1F(d) in~\cite{Engelking} and Theorem 4.2 in~\cite{Hodel}). The example given in~\ref{kappaO}, where $X=\kappa\omega$, is a separable Hausdorff space with countable pseudocharacter. Thus, $|X|=2^\mathfrak{c}>\mathfrak{c}=2^{d(X)\psi(X)}$, showing $2^{d(X)\psi(X)}$ is not a bound for the cardinality of all Hausdorff spaces. This is in contrast to Theorem~\ref{dwpsi}, which states that $2^{d(X)w\psi_c(X)}$ (and $2^{d(X)\psi_c(X)}$) is in fact a bound for the cardinality of all Hausdorff spaces.

The previous theorem has an immediate and well-known corollary.

\begin{corollary} If $X$ is Hausdorff then $|X|\leq 2^{2^{d(X)}}$. 
\end{corollary}

\begin{proof}
By Theorem~\ref{dwpsi} and Proposition~\ref{wpc3}, we have
$$|X|\leq 2^{d(X)w\psi_c(X)}\leq 2^{d(X)\cdot 2^{d(X)}}=2^{2^{d(X)}}.$$
\end{proof}

We arrive now at two new cardinality bounds for Hausdorff spaces that do not involve any variation of the pseudocharacter.

\begin{corollary}\label{piwdot}
Let $X$ be a Hausdorff space. Then, 
\begin{enumerate}
\item $|X|\leq 2^{\pi w(X)^{dot(X)}}$, and 
\item $|X|\leq 2^{d(X)^{\pi\chi(X)}}$.
\end{enumerate}
\end{corollary}

\begin{proof}
For (1), by Theorem~\ref{dwpsi} and Proposition~\ref{theProp}(8), we have
$$|X|\leq 2^{d(X)w\psi_c(X)}\leq 2^{d(X)\pi w(X)^{dot(X)}}=2^{\pi w(X)^{dot(X)}}.$$
For (2), by (1) and the fact that $dot(X)\leq\pi\chi(X)$ for any space, we have
$$|X|\leq 2^{\pi w(X)^{dot(X)}}=2^{(d(X)\pi\chi(X))^{dot(X)}}\leq 2^{{d(X)\pi\chi(X)}^{\pi\chi(X)}}=2^{d(X)^{\pi\chi(X)}}.$$
\end{proof}

Note that the placement of $dot(X)$ as an exponent in Corollary~\ref{piwdot}(1) is necessary. For example, consider $X=\beta\omega$, the Stone-\v Cech compactification of the natural numbers. $X$ has a countable dense set of isolated points and so $\pi w(X)=\omega$ and $dot(X)\leq\pi\chi(X)\leq\pi w(X)=\omega$. Therefore $|X|= 2^\mathfrak{c}>\mathfrak{c}=2^{\pi w(X)dot(X)}$, indicating that dropping the placement of $dot(X)$ is invalid in Corollary~\ref{piwdot}(1). Nonetheless, one should consider $dot(X)$ a ``small" cardinal function as $dot(X)\leq\min\{ot(X),\pi\chi(X)\}$ by Proposition~\ref{theProp}(1) and (2).

In~\cite{GTT16} it was shown by Gotchev, Tkachenko, and Tkachuk that if $X$ is Hausdorff then $|X|\leq\pi w(X)^{ot(X)\psi_c(X)}$. We show in the next result that $\psi_c(X)$ can be replaced with $w\psi_c(X)$.

\begin{theorem}\label{pwot}
If $X$ is Hausdorff then $|X|\leq\pi w(X)^{ot(X)w\psi_c(X)}$.
\end{theorem}

\begin{proof}
Let $\scr{B}$ be a $\pi$-base such that $|\scr{B}|=\pi w(X)$ and let $\kappa=ot(X)w\psi_c(X)$. For all $x\in X$ let $\scr{V}_x$ be a weak closed pseudobase at $x$ such that $|\scr{V}_x|\leq\kappa$. For all $x\in X$ and $V\in\scr{V}_x$, let $\scr{B}(x,V)=\{B\in\scr{B}:B\sse V\}$. Then $x\in\cl{V}=\cl{\Un\scr{B}(x,V)}$. As $ot(X)\leq\kappa$, there exists $\scr{C}(x,V)\in[\scr{B}(x,V)]^{\leq\kappa}$ such that $x\in\cl{\Un\scr{C}(x,V)}$. (Note that in this situation we cannot use Proposition~\ref{dot} as $x$ may not be in $V$). Then $\{x\}=\Meet_{V\in\scr{V}_x}\cl{\Un\scr{C}(x,V)}$. Enumerate $\scr{V}_x=\{V(x,\alpha):\alpha<\kappa\}$.

Define $\phi:X\to\left([\scr{B}]^{\leq\kappa}\right)^{\kappa}$ by $\phi(x)(\alpha)=\scr{C}(x,V(x,\alpha))$. We show $\phi$ is one-to-one. Suppose $x$ and $y$ are two distinct points in $X$. Then 
$$\Meet_{\alpha<\kappa}\cl{\Un\scr{C}(x,V(x,\alpha))}\neq\Meet_{\alpha<\kappa}\cl{\Un\scr{C}(y,V(y,\alpha))}.$$
It follows that there exists $\alpha<\kappa$ such that $\scr{C}(x,V(x,\alpha))\neq\scr{C}(y,V(y,\alpha))$ and thus $\phi(x)(\alpha)\neq\phi(y)(\alpha)$. This shows $\phi(x)\neq\phi(y)$ and that $\phi$ is one-to-one. It follows that 
$$|X|\leq\left|\left([\scr{B}]^{\leq\kappa}\right)^{\kappa}\right|\leq\pi w(X)^\kappa=\pi w(X)^{ot(X)w\psi_c(X)}.$$
\end{proof}

Using the above Theorem~\ref{pwot} we obtain another new bound for the cardinality of any Hausdorff space that does not involve any notion of pseudocharacter.

\begin{corollary}\label{bound42}
If $X$ is Hausdorff then $|X|\leq d(X)^{\pi\chi(X)^{ot(X)}}$.
\end{corollary}

\begin{proof}
By Theorem~\ref{pwot} and Proposition~\ref{wpc1}, we have
\begin{align}
|X|&\leq\pi w(X)^{ot(X)w\psi_c(X)}\leq\pi w(X)^{ot(X)\pi\chi(X)^{ot(X)}}=\pi w(X)^{\pi\chi(X)^{ot(X)}}\notag\\
&=(d(X)\pi\chi(X))^{\pi\chi(X)^{ot(X)}}=d(X)^{\pi\chi(X)^{ot(X)}}.\notag
\end{align}
\end{proof}

Pospi{\v s}il~\cite{Pos37} showed in 1937 that $|X|\leq d(X)^{\chi(X)}$ for any Hausdorff space. Bella and Cammaroto~\cite{BelCam88} improved this bound to $|X|\leq d(X)^{t(X)\psi_c(X)}$ in 1988. In 2018 the author showed in~\cite{Car2018} that $t(X)$ can be replaced with $wt(X)$. We show in Theorem~\ref{densityWt} that in fact $\psi_c(X)$ can further be replaced with $w\psi_c(X)$. First we need a definition and some lemmas.

In \cite{JVM2018}, \juhasz~and van Mill introduced the notion of a $\scr{C}$-saturated subset of a space $X$. 
\begin{definition}
Given a cover $\scr{C}$ of $X$, a subset $A\sse X$ is $\scr{C}$-\emph{saturated} if $A\meet C$ is dense in $A$ for every $C\in\scr{C}$. 
\end{definition}
It is clear that the union of $\scr{C}$-saturated subsets is $\scr{C}$-saturated. The following was given in~\cite{JVM2018} in the case $\kappa=\omega$, and extended to the general case in~\cite{Car2018}.

\begin{lemma}[\cite{JVM2018},\cite{Car2018}]\label{saturated}
Let $X$ be a space, $wt(X)=\kappa$, and let $\scr{C}$ be a cover witnessing that $wt(X)=\kappa$. Then for all $x\in X$ there exists $S_x\in[X]^{\leq 2^\kappa}$ such that $x\in S_x$ and $S_x$ is $\scr{C}$-saturated.
\end{lemma}

\begin{lemma}[Proposition 2.4 in~\cite{Car2018}]\label{cardlemma}
Let $X$ be a space, $D\sse X$, and suppose there exists a cardinal $\kappa$ such that for all $x\in X$ there exists $\scr{B}_x\in\left[[D]^{\leq\kappa}\right]^{\leq\kappa}$ such that $\{x\}=\Meet_{B\in\scr{B}_x}\cl{B}$. Then $|X|\leq |D|^{\kappa}$.
\end{lemma}

\begin{theorem}\label{densityWt}
If $X$ is Hausdorff then $|X|\leq d(X)^{wt(X)w\psi_c(X)}$.
\end{theorem}

\begin{proof} 
Let $\kappa=wt(X)w\psi_c(X)$ and let $D$ be a dense subset such that $|D|=d(X)$. For all $x\in X$ let $\scr{V}$ be a weak closed pseudobase such that $|\scr{V}_x|\leq\kappa$. Let $\scr{C}$ be a cover of $X$ witnessing that $wt(X)\leq\kappa$. By Lemma~\ref{saturated}, for all $x\in X$ there exists a $\scr{C}$-saturated set $S_x$ such that $x\in S_x$ and $|S_x|\leq 2^\kappa$. Let $S=\Un_{d\in D}S_d$. Then $S$ is $\scr{C}$-saturated, as $S$ is the union of $\scr{C}$-saturated sets, and $|S|\leq |D|\cdot 2^\kappa=|D|^\kappa$. Observe that as $D\sse S$, we have that $S$ is dense in $X$.

Fix $x\in X$ and let $C\in\scr{C}$ such that $x\in C$. We show for each $V\in\scr{V}_x$ that $x\in cl_C(V\meet S\meet C)$. As $S$ is $\scr{C}$-saturated and dense in $X$, we have
$$x\in\cl{V}=\cl{V\meet S}=cl(V\meet cl_S(S\meet C))\sse cl(V\meet cl(S\meet C))=cl(V\meet S\meet C).$$

Therefore $x\in C\meet cl(V\meet S\meet C)=cl_C(V\meet S\meet C)$. As $t(C)\leq\kappa$, there exists $A_V\sse V\meet S\meet C$ such that $x\in cl_C(A_V)\sse cl(A_V)$ and $|A_V|\leq\kappa$.

Then,
$$\{x\}\sse\Meet_{V\in\scr{V}_x}cl(A_V)\sse\Meet_{V\in\scr{V}_x}cl(V\meet S\meet C)\sse\Meet_{V\in\scr{V}_x}\cl{V}=\{x\}.$$

This shows $\{x\}=\Meet_{V\in\scr{V}_x}cl(A_V)$. Now, observe that $\scr{B}_x=\{A_V:V\in\scr{V}_x\}\in\left[[S]^{\leq\kappa}\right]^{\leq\kappa}$. By Lemma~\ref{cardlemma}, we have $|X|\leq|S|^\kappa\leq\left(|D|^\kappa\right)^\kappa=|D|^\kappa=d(X)^{wt(X)w\psi_c(X)}$.

\end{proof}

\section{Remarks concerning homogeneous spaces.}

Recall that a space $X$ is \emph{homogeneous} if for all $x,y\in X$ there exists a homeomorphism $h:X\to X$ such that $h(x)=y$ (Definition~\ref{definehomogeneous}). 

Using that $|X|\leq 2^{\psi(X)}$ for any compact Hausdorff space $X$, it can easily be seen that the cardinality of a $\sigma$-compact space $X$ is bounded above by $\omega\cdot 2^{\psi(X)}=2^{\psi(X)}$. The next example shows that $2^{w\psi_c(X)}$ is not a bound for the cardinality of all $\sigma$-compact homogeneous spaces $X$.

\begin{example}\label{homogeneousexample}
In~\cite{CR2012}, the author and Ridderbos demonstrated that for every infinite cardinal $\kappa$ there is a $\sigma$-compact homogeneous space of cardinality $\kappa$ with countable tightness and countable $\pi$-character. (Note that if such spaces were compact, then their cardinality would be bounded by $\mathfrak{c}$. This is due to de la Vega's result that the cardinality of a compact homogeneous space is at most $2^{t(X)}$~\cite{DLV2006}). We refer the reader to the very lengthy description of these examples in~\cite{CR2012}. 

For each infinite cardinal $\kappa$, let $X_\kappa$ be a $\sigma$-compact homogeneous space of cardinality $\kappa$ with countable tightness and countable $\pi$-character. We observe that, by Proposition~\ref{wpc1}, for each $X_\kappa$ we have $w\psi_c(X_\kappa)\leq\pi\chi(X_\kappa)^{ot(X_\kappa)}\leq\omega^\omega=\mathfrak{c}$. Since each $X_\kappa$ is $\sigma$-compact, it is the countable union of compact subspaces. Thus, $|X_\kappa|\leq\omega\cdot 2^{\psi(X_\kappa)}=2^{\psi(X_\kappa)}$. However, when $\lambda=\left(2^\mathfrak{c}\right)^+$, we have $\left|X_\lambda\right|=\left(2^\mathfrak{c}\right)^+>2^{\mathfrak{c}}\geq 2^{w\psi_c(X_\lambda)}$. Therefore, $2^{w\psi_c(X)}$ is not a bound for the cardinality of all $\sigma$-compact homogeneous spaces. Additionally, as $\left(2^\mathfrak{c}\right)^+=|X_\lambda|\leq 2^{\psi(X_\lambda)}$, we see that $\psi(X_\lambda)>\mathfrak{c}$ while $w\psi_c(X_\lambda)\leq\mathfrak{c}$. In this manner by choosing $\kappa$ large enough we can make $\psi(X_\kappa)$ as large as we want while $w\psi_c(X_\kappa)\leq\mathfrak{c}$. This shows that the spread between $w\psi_c(X)$ and $\psi(X)$ can be arbitrarily large among $\sigma$-compact homogeneous spaces. 
\end{example}

However, even given Examples~\ref{compactexample} and~\ref{homogeneousexample} it is still unclear whether $2^{w\psi_c(X)}$ is a bound for the cardinality of all compact, homogeneous Hausdorff spaces. (The examples given in~\ref{compactexample} are compact but not homogeneous, and the examples given in~\ref{homogeneousexample} are homogeneous but not compact). So we ask the following.

\begin{question}
If $X$ is a compact, homogeneous, Hausdorff space, is $|X|\leq 2^{w\psi_c(X)}$?
\end{question}

As $w\psi_c(X)\leq\pi\chi(X)^{ot(X)}$ by Proposition~\ref{wpc1}, we can also ask the following:

\begin{question}
If $X$ is a compact, homogeneous, Hausdorff space, is $|X|\leq 2^{\pi\chi(X)^{ot(X)}}$?
\end{question}

A stronger version of this question is as follows.

\begin{question}
If $X$ is a compact, homogeneous, Hausdorff space, is $|X|\leq 2^{ot(X)\pi\chi(X)}$?
\end{question}

A positive answer to the above question would simultaneously improve the result that $|X|\leq 2^{wt(X)\pi\chi(X)}$ for homogeneous compact Hausdorff spaces (see~\cite{BC2020}) and the result that $|X|\leq 2^{c(X)\pi\chi(X)}$ for any homogeneous Hausdorff space (\cite{CR2008}), as $ot(X)\leq\min\{wt(X),c(X)\}$.

The invariant $q\psi(X)$ (Definition~\ref{qpsi}) was defined by Ismail in~\cite{Ism81} for a Hausdorff space $X$. Ismail showed the following.

\begin{proposition}[\cite{Ism81}]\label{IsHomog} If $X$ is a homogeneous Hausdorff space then $|X|\leq |RO(X)|^{q\psi(X)}$.
\end{proposition}

Notice that this improves the result that if $X$ is Hausdorff then $|X|\leq |RO(X)|^{w\psi_c(X)}$ (Proposition~\ref{ROcardBound}) in the homogeneous case, as $q\psi(X)\leq w\psi_c(X)$ for any Hausdorff space by Proposition~\ref{theProp}(6).

By Corollary~\ref{QRRObound}, we have $|RO(X)|\leq\pi\chi(X)^{c(X)}$ for a quasiregular Hausdorff space. Combining that with Proposition~\ref{IsHomog}, we have the following result. This was first mentioned in~\cite{Car2021} for regular homogeneous spaces but in fact works for quasiregular homogeneous spaces.

\begin{theorem}
If $X$ is a homogeneous quasiregular Hausdorff space then $|X|\leq\pi\chi(X)^{c(X)q\psi(X)}$.
\end{theorem}

Observe that this improves the cardinality bound $2^{c(X)\pi\chi(X)}$ for homogeneous Hausdorff spaces given in~\cite{CR2008} in the case when $X$ is additionally quasiregular, as $q\psi(X)\leq\pi\chi(X)$ by Proposition~\ref{theProp}(7). 

Using the fact that $|RO(X)|\leq\pi\chi(X)^{c(X)d\psi_c(X)}$ (Corollary~\ref{HausRO}) along with~\ref{IsHomog}, we obtain the result that if $X$ is a homogeneous Hausdorff space then $|X|\leq\pi\chi(X)^{c(X)d\psi_c(X)q\psi(X)}$. This would seem to be a variation of the bound $2^{c(X)\pi\chi(X)}$. Yet it is easily seen that for a homogeneous space $X$ if there is a point with closed pseudocharacter $\kappa$, then $\psi_c(X)\leq\kappa$. Thus if $X$ is homogeneous and Hausdorff we have $d\psi_c(X)=\psi_c(X)$. Therefore the result is equivalent to this: If $X$ is homogeneous and Hausdorff then $|X|\leq\pi\chi(X)^{c(X)\psi_c(X)q\psi(X)}=\pi\chi(X)^{c(X)\psi_c(X)}$. However, this is exactly Sun's inequality $|X|\leq\pi\chi(X)^{c(X)\psi_c(X)}$ for any Hausdorff space $X$. Therefore the bound $2^{c(X)\pi\chi(X)}$~\cite{CR2008} for the cardinality of a homogeneous Hausdorff space $X$ seems to be the strongest bound in this direction.

Combining Ismail's result in Proposition~\ref{IsHomog} with the two bounds for $|RO(X)|$ given in Proposition~\ref{theProp}(8), we obtain two new bounds for the cardinality of a homogeneous Hausdorff space.

\begin{theorem}\label{newhomog} Let $X$ be a homogeneous Hausdorff space. Then,
\begin{itemize}
\item[(1)] $|X|\leq 2^{d(X)q\psi(X)}$, and
\item[(2)] $|X|\leq \pi w(X)^{c(X)q\psi(X)}$.
\end{itemize}
\end{theorem}

It was shown independently in~\cite{DLVthesis} and~\cite{RidMasters} that if $X$ is a homogeneous Hausdorff space then $|X|\leq d(X)^{\pi\chi(X)}$. (Also, see~\cite{Car2021}). Observe that this bound also does not contain $\psi_c(X)$ or its variations. Compare this bound with Theorem~\ref{piwdot}(2), which states that the cardinality of any Hausdorff space is at most $2^{d(X)^{\pi\chi(X)}}$. Additionally compare this with Corollary~\ref{bound42}, which gives $d(X)^{\pi\chi(X)^{ot(X)}}$ as a bound for the cardinality of any Hausdorff space, and Theorem~\ref{newhomog}(1). Finally, the bound $2^{c(X)\pi\chi(X)}$ for the cardinality of homogeneous Hausdorff space should be compared to the bounds $2^{c(X)\pi\chi(X)^{ot(X)}}$ and $2^{c(X)^{\pi\chi(X)}}$ for quasiregular Hausdorff spaces given in Theorem~\ref{QRBound}.


\begin{thebibliography}{21}

\bibitem{BelCam88}
A.~Bella and F.~Cammaroto, \emph{On  the cardinality of {U}rysohn
spaces}, 
Canad. Math. Bull. \textbf{31} (1988), no.~2, 153--158.

\bibitem{BC2020}A.~ Bella, N.~ Carlson, \emph{On weakening tightness to weak tightness}, 
Monatsh. Math. \textbf{192}, no. 1 (2020), 39--48.

\bibitem{BCG2022a} A.~Bella, N.~Carlson, I.~Gotchev, \emph{More on cardinality bounds involving the weak Lindel\"of degree}, Quaestiones Mathematicae, DOI: 10.2989/16073606.2022.2040634

\bibitem{BCG2022b} 
A.~Bella, N.~Carlson, I.~Gotchev, \emph{On spaces with a $\pi$-base with elements with compact closure}, preprint.

\bibitem{Car2007}
N.~Carlson, \emph{Non-regular power homogeneous spaces}, Topology Appl.
  \textbf{154} (2007), no.~2, 302--308.
  
\bibitem{Car2018} N.~Carlson, \emph{On the weak tightness and power homogeneous compacta}, 
Topology Appl. \textbf{249} (2018), 103--111.

\bibitem{Car2021} N.~ Carlson, \emph{A survey of cardinality bounds on homogeneous topological spaces}, Top. Proc. \textbf{57} (2021), 259--278.

\bibitem{CR2008}
N.~Carlson and G.~J. Ridderbos, \emph{Partition relations and power homogeneity}, Top. Proc. \textbf{32} (2008), 115--124.

\bibitem{CR2012}N.~Carlson, G.J.~Ridderbos, \emph{On several cardinality bounds on power homogeneous spaces}, Houston Journal of Mathematics, \textbf{38} (2012), no.~1, 311-332.

\bibitem{Cha77}
A.~Charlesworth,  \emph{On  the   cardinality  of  a  topological space}, 
Proc. Amer. Math. Soc. \textbf{66} (1977), 138--142.

\bibitem{Efimov} B.~Efimov, \emph{Extremally disconnected bicompacta of continuum $\pi$-weight}, Soviet Math. Dokl. 9
(1968), 1404-1407.

\bibitem{Engelking}
R. Engelking, \emph{General {T}opology}, Heldermann Verlag, Berlin, second ed., 1989.

\bibitem{GTT16} I.~Gotchev, M.~Tkachenko, V.~Tkachuk, 
\emph{Regular $G_\delta$-diagonals and some upper bounds for cardinality of topological spaces}, 
Acta Math. Hungar. \textbf{149} (2016), no. 2, 324–337. 

\bibitem{DeG1965} J.~De Groot, \emph{Discrete subspaces of Hausdorff spaces}, Bull. Acad. Polon. Sci. \textbf{13} (1965), 537--544.

\bibitem{Hodel}
R.~Hodel, \emph{Cardinal functions I}, in Handbook of Set-theoretic Topology (K. Kunen and J.E. Vaughan, eds), North-Holland, Amsterdam, 1984, 1--61.

\bibitem{Ism81}
M.~Ismail, \emph{Cardinal functions of homogeneous spaces and topological groups}, Math.
Japon. 26 (1981), no. 6, 635?646.

\bibitem{Juh80}
I.~Juh{\'a}sz, \emph{Cardinal {F}unctions in {T}opology---{T}en
{Y}ears
  {L}ater}, second ed., Mathematical Centre Tracts, vol. 123,
Mathematisch
  Centrum, Amsterdam, 1980.
  
\bibitem{JVM2018}
I.~\juhasz, J.~van Mill, \emph{On $\sigma$-countably tight
spaces}, Proc. Amer. Math. Soc. \textbf{146} (2018), no. 1,
429--437.

\bibitem{Pos37}
B.~Pospi{\v s}il, \emph{Sur la puissance d'un espace contenant
une partie dense
  de puissance donn\'ee}, \v{C}asopis Pestov\'ani Matematiky a
Fysiky \textbf{67} (1937), 89--96.

\bibitem{RidMasters}
G.~J. Ridderbos, \emph{Notes on cardinal functions and homogeneity}, Master Thesis, Vrije Universiteit, Amsterdam (unpublished).

\bibitem{Sap1972}
B.~\sapirovskii, \emph{On discrete subspaces of topological spaces; weight, tightness and Suslin number}, Soviet
Math. Dokl. 13 (1972) 215--219.

\bibitem{Sap1974}
B.~\sapirovskii, \emph{Canonical sets and character. Density and weight in bicompacta} (Russian), Dokl. Akad. Nauk SSSR 218 (1974), 58–61.

\bibitem{Sun88}S.~H.~Sun, \emph{Two new topological cardinal inequalities}, 
Proc. Amer. Math. Soc. \textbf{104} (1988), 313--316.

\bibitem{Tka83}
M.\,G.~Tka\v{c}enko,   \emph{The notion of o-tightness and C-embedded subspaces of  products}, 
Topology  Appl.  \textbf{15} (1983), 93--98.

\bibitem{DLVthesis}
R.~de~la Vega, \emph{Homogeneity properties on compact spaces}, Doctoral Thesis, University of Wisconsin-Madison.

\bibitem{DLV2006}
R.~de~la Vega, \emph{A new bound on the cardinality of homogeneous compacta}, Topology Appl. \textbf{153} (2006), 2118-2123.

\end{thebibliography}
\end{document}